\documentclass[a4paper,10pt]{article}

\def \Z {{\mathbf {Z}}}
\def \R {{\mathbf {R}}}

\title{ Local rank of ergodic symmetric $n$-powers does not exceed $n!n^{-n}$}
\author{V.V. Ryzhikov}
\date{vryzh@mail.ru}

\begin{document}
\Large
\maketitle

\section{Introduction}

Let $S$ be an automorphism of  a measure space $(X,\mu)$, $\mu(X)=1$. We denote by  
 $\beta(S)$ its \it  local rank. \rm  It is defined as   maximal   $\beta$ such that
there is a partition sequence 
$$ \xi_j=\{ B_j,  SB_j,    S^2B_j,    \dots,  S^{h_j-1}B_j,C_j^1,\dots, C_j^{m_j}\dots\}$$
for which 
any measurable set can be approximated by 
$\xi_j$-measurable ones as $j\to\infty$, and 
 $\mu(U_j)\to \beta$, where 
$U_j=\bigsqcup_{0\leq k<h_j}S^kB_j$ are  called $\beta$-towers.

{\bf THEOREM.}\footnote
{Our results have been presented at  the conference
MODERN THEORY OF DYNAMICAL SYSTEMS AND APPLICATIONS TO THEORETICAL CELESTIAL MECHANICS dedicated to the memory and the 70th birthday of V.M. Alexeyev, Moscow, December 2002. } \it   Let  $T$  be weakly mixing, then

\bf a) $\beta(T^{\otimes n})\leq n^{-n}$, \ \  \bf b) $\beta(T^{\odot n})\leq n!n^{-n}$,
 \ 
\
c) $ Rank(T^{\odot n})=\infty$ as $n>1$. \rm
\vspace{2mm}

A.Katok \cite{K} showed that for a generic $T$ it holds
 $\beta(T^{\times n})\geq n^{-n}$,
 $\beta(T^{\odot n})\geq n!n^{-n}$, so, the above bounds are exact.  
In fact the above theorem is true for all $T$, but we shall consider  
only weakly mixing automorphisms,   generalizing  \cite{R1}.
\section{Auxiliary assertions}

{\bf LEMMA 1. \cite{R1} } \it Let an ergodic  automorphism $S$
of a space   $(\bar{X},\bar{\mu})$ commute with  an automorphism
$R$, and    $\beta(S)>0$.
Then for any   $\delta >0$
there is  $m>0$  such that  the following  weak convergence 
$$  \hat{Y}_j\circ S^{n_{j}}\to  (1-\delta')\beta(S) {R}^m
$$
holds for   some  sequence $n_{j}$ and    operators   $  \hat{Y}_j$ of multiplication by certain  sets  ${Y_j}$.   Here  ${Y_j}$ are sub-towers of our $\beta$-towers, 
 $\bar{\mu}({Y}_j)\to (1-\delta')\beta(S)$ for some  $\delta'\in [0,\delta]$.
\rm
\medskip

{ Proof.}
Let  $U_j$ be a sequence of $\beta$-towers for  $S$.
Let's define small towers
$ U_j^{\delta}=\bigsqcup_{0\leq k \leq \delta h_j}S^kB_j$.
For some  $m>0$  we have
$$     \limsup_{j} \bar{\mu}(S^m U_{j}^{\delta} \cap  U_j^{\delta}) > 0.$$
Indeed, from  $ \mu(U_j^{\delta})\approx {\delta}\beta(S)$ we get that  the sets
$U_j^{\delta},  R U_j^{\delta},\dots,   R^N U_j^{\delta}$
are not disjoint as   ${\delta}\beta(S) N > 1$.
So, for some  $m$, $0<m\leq N$, and   $c>0$   the inequality 
$$\bar{\mu} (U_{j}^{\delta} \cap R^m U_{j}^{\delta}) >c $$
 holds for an infinite collection  of $j$-s.
From \cite{R},  \S 3,  it follows that an ergodic (!)  joining
$\Delta_{R^m}$ corresponding to the operator 
 $R^m$ can be approximated by  parts of off-diagonal measures $\Delta_{S^{n_j}}$, and these parts are situated in
 $ Y_j\times \bar{X}$, where
$$Y_j=\bigsqcup_{\delta' h_j\leq k \leq h_j}S^kB_j$$ for 
$\delta'\leq\delta$.
The latter is equivalent to the assertion of our lemma.     
\medskip

{\bf LEMMA 2.}   \it If  $\delta>0$ and   
 $\beta(T^{\otimes n})>0$, then for some  $Y_j$ (sub-towers of our $\beta$-towers)   there is a sequence  $k_j\to\infty$ such that 
for $\delta'<\delta$
$$ \hat{Y}_j\circ (T^{\otimes n})^{k_{j}}\to (\beta(T^{\otimes n})-\delta')( T^m \otimes T^{2m}\otimes\dots \otimes T^{nm}).$$

\medskip                \rm

Proof. We apply Lemma 1 for 
$S= T^{\otimes n}$, $R= T \otimes T^{2}\otimes\dots \otimes T^{n}$.

Further we use the following notation:
$$P_1\odot\dots\odot P_n =\frac{1}{n!}\sum_{\sigma}P_{\sigma(1)}\otimes\dots\otimes P_{\sigma(n)},$$ 
where $\sigma$ runs all permutations.

{\bf LEMMA 3.}   \it If  $\delta>0$ and  
 $\beta(T^{\odot n})>0$, then for a sequence  $Y_j$ and some  $\delta'<\delta$ there is a sequence  $k_j\to\infty$ such that 
$$ \hat{Y}_j\circ (T^{\odot n})^{k_{j}}\to (\beta(T^{\odot n})-\delta')( T^m \odot T^{2m}\odot\dots \odot T^{nm}).$$\rm

Proof.  We consider on $X^{\odot n}\times X^{\odot n}$  a joining $\nu$  (a finite-value polymorphism in Vershik's terminology \cite{V})  corresponding  to the 
operator  $ P_m=T^m \odot T^{2m}\odot\dots \odot T^{nm}$.  It is
ergodic. Indeed, it is clear that $\nu$ as a joining on $X^n \times X^n$  is ergodic with respect to the tensor square of the  action
of $T^{\otimes n}$ and  all coordinate permutations.   But  the permutations act on symmetric sets  identically.  So, our joining is ergodic with respect to $T^{\odot n}\times T^{\odot n}$, and we may act as in  Lemma 1 taking into account  the following remark.  We find $m$ for which
$$\bar{\mu} (U_{j}^{\delta} \cap R^m U_{j}^{\delta}) >c'>0, $$
thus, we get
$${\nu} (U_{j}^{\delta} \times U_{j}^{\delta}) =
\int \chi_{U_{j}^{\delta}} P_m\chi_{U_{j}^{\delta}}\  d \mu^{n} \ >c>0. $$
A nature of a joining is not  important, but its ergodicity is essential. The mentioned approximations from  \cite{R},  \S 3, lead us to the assertion of Lemma 3.

\section{ Upper bound  for local rank  of  symmetric powers}
We prove a) and b).
Let (from Lemma 2)
$$ \hat{Y}_j\circ (T^{\otimes n})^{k_{j}}\to (\beta(T^{\otimes n})-\delta')( T^m \otimes T^{2m}\otimes\dots \otimes T^{nm}),\eqno (1)$$
and (by a choice of subsequence)
$$T^{k_j}\to aQ+\sum_{p=1}^n a_pT^{mp},\eqno (2)$$ 
where the Markov operator $Q\perp  T^{mp}$, $m=1,2,...n$ (the corresponding
polymorphisms are disjoint as measures on  $X^n\times X^n$).
From (2) we have 
$$(T^{\otimes n})^{k_j}\to \left(aQ+\sum_{p=1}^n a_pT^{mp}\right)^{\otimes n}=a_1a_2\dots a_n T^m \otimes T^{2m}\otimes\dots \otimes T^{nm}+ b P,\eqno (3)$$ 
where the Markov operator $P\perp P_m=T^m \otimes T^{2m}\otimes\dots \otimes T^{nm}$ (an exercise). 
Thus, comparing (3) with (1), we get
$$\beta(T^{\otimes n})-\delta'\leq a_1a_2\dots a_n \leq n^{-n}$$
(we note that $a_1+a_2+\dots+a_n\leq 1$).
Since $\delta'$ is arbitrary small, we obtain
$$\beta(T^{\otimes n}) \leq n^{-n}.$$

So, a) is proved. Now we prove b). Let
$$ \hat{Y}_j\circ (T^{\odot n})^{k_{j}}\to (\beta(T^{\odot n})-\delta')( T^m \odot T^{2m}\odot\dots \odot T^{nm}),$$
and (by a choice of subsequence)
$$T^{k_j}\to aQ+\sum_{p=1}^n a_pT^{mp}, \ m>0, \eqno (4)$$ 
where the Markov operator $Q\perp  T^{mp}$, $p=1,2,...n$. 
We have
$$(T^{\otimes n})^{k_j}\to \left(aQ+\sum_{p=1}^n a_pT^{mp}\right)^{\otimes n}=n!a_1a_2\dots a_n T^m \odot T^{2m}\odot\dots \odot T^{nm}+ b P.$$ 
If $P$ and $P_m=T^m \odot T^{2m}\odot\dots \odot T^{nm}$ both restricted to $L_2(X^{\odot n})$  are not disjoint (in the above sense), then extended 
$P$ (as an operator in $L_2(X^{ n})$ )  is simply obliged to have  $P_m$ as a component . 
Indeed, let  $P=bP'+\dots$, where $P'$ commutes with all coordinate permutations and the restriction of $P'$ to $L_2(X^{\odot n})$ coincides with $P_m$.
 The corresponding joining $\nu'$ is a self-joining for the mentioned group
action,  it is absolutely continuous with respect to the ergodic self-joining in $X^n\times X^n$ that corresponds to $P_m$.
Thus, $\nu'=\nu$, $P'=P_m$.  

The operator $P$ is a convex sum of products containing $Q$ as a multiplier, say 
$  T^m\otimes Q\otimes\dots$.  If $P$ possesses  the component
$T^m \otimes T^{2m}\otimes\dots \otimes T^{nm}$, then
$Q$  have a component $T^{2m}$. This  contradicts the definition of $Q$ from (4).  

Thus, $\beta(T^{\odot n})-\delta'\leq n!a_1a_2\dots a_n$,  hence,  
 $\beta(T^{\odot n}) \leq n!n^{-n}.$

\section{Infinity of Rank for symmetric powers }
If
$$T^{k_j}\to aQ+\sum_{p=1}^n a_pT^{mp},\ m>0,$$ 
and the  operator $Q\perp  T^{mp}$, $p=1,2,...n$,
then there is a collection of sets $X^m_j$, $p=1,2,...n$, such that
$$ \hat X^m_j T^{k_j}\to a_pT^{pm}, \ j\to\infty$$
(a proof based on  a  simple  technique of polymorphisms is given in \cite{R2}, Theorem 7.1).

If in addition 
$$ \hat{Y}_j\circ (T^{\otimes n})^{k_{j}}\to \beta T^m \otimes T^{2m}\otimes\dots \otimes T^{nm}$$
holds, then we get (as an obvious consequence)
$$\mu^n(Y_j\setminus (X^1_j\times\dots\times X^n_j))\to 0.$$
(From this and Lemma 2 we prove  again a).)

 Our aim is to prove the infinity of usual Rank, i.e. to show that  \it for any $r$
 the space $X^{\odot n}$ cannot be asymptotically covered with a vanishing error  by a  collection
of $\beta_i$-towers  $U_j^i$, $i=1, 2, \dots , r$. \rm

Let's consider $\beta_i$-towers  as symmetric sets on 
$X^n$.
If for  symmetric sets $Y_j$ 
$$ \hat{Y}_j\circ (T^{\odot n})^{k_{j}}\to \beta T^m \odot T^{2m}\odot\dots \odot T^{nm},$$
then we get
$$\mu^n\left(Y_j\,\setminus \bigcup_\sigma (X^{\sigma(1)}_j\times\dots\times X^{\sigma(n)}_j)\right)\to 0.$$
From this fact and  Lemma 3 we deduce that $\beta_i$-towers $U_{j'}^i$
are asymptotically  situated (with a vanishing error) in  sets 
$Z_{j'}=\bigcup_\sigma (X^{\sigma(1)}_{j'}\times\dots\times X^{\sigma(n)}_{j'})$, where    $X^p_{j'}$, $p=1,2,...,n$, \bf are
pairwise disjoint. \rm It is not hard to see that  we can't cover (with a vanishing error)
the space $X^n$  by a union of a finite collection of such special symmetric sets $Z_j$, i.e  $T^{\odot n}$ is not of finite Rank.
\vspace{2mm}

\bf Remark. \it Our statements (on $\Z$-actions)  are also true for $\Z^d$-actions (as well as for $\R^d$-flows). We note also that
 for any $d\geq 1$ and $b\in (0,n!n^{-n})$ there is a partially mixing rank one action of $\Z^d$ ($\R^d$)  such that $\beta({\bf action}^{\odot n})=b.$ \rm

 \end{document}